# General Divisibility Criteria

Anatoly A. Grinberg[1] and Serge Luryi

*Electrical and Computer Engineering Department*
*University at Stony Brook, Stony Brook, NY 11794-2350*

ABSTRACT

We believe we have made progress in the age-old problem of divisibility rules for integers. Universal divisibility rule is introduced for any divisor in any base number system. The divisibility criterion is written down explicitly as a linear form in the digits of the test number. The universal criterion contains only two parameters that depend on the divisor and are easily calculated. These divisibility parameters are not unique for a given divisor but each two-parameter set yields a unique criterion. Well-known divisibility rules for exemplary divisors in the decimal system follow from the universal expression as special cases.

**I. Introduction**

Divisibility rules are designed to answer the question of the divisibility of a test integer $A$ by a divisor $n$ — without actually performing the division. The usual rule corresponds to forming a criterion number $C$ that is smaller than the test number (ideally, $C << A$) and possesses the same property in terms of the divisibility by $n$.

Divisibility rules have been derived for many divisors [1, 2]. Some of these rules involve only a few rightmost digits of the test number. Thus the well-known criteria for divisibility into *2, 4,* or *8* involve (in the decimal system) the *1, 2,* or *3* rightmost digits, respectively.

We shall be concerned with the rules involving *all* digits of the test number and presented as a linear form in these digits. For an (*m+1*)-digit test number *A*,

$$A = \sum_{k=0}^{m} t^k a_k , \qquad (I.1)$$

where $a_k$ are its digits in the *t*-base number system, we seek the criterion of divisibility of *A* into *n* in the form

$$C = \sum_{k=0}^{m} c_k(n) a_k , \qquad (I.2)$$

---

[1] Corresponding author: anatoly_gr@yahoo.com

Jan 21, 2014

where $c_k$ are the coefficients determined in terms of the divisor. Examples of such rules in the decimal system ($t=10$) are the criteria for $n = 3$ or $11$. Thus, the commonly familiar rule for division into *3* corresponds to all $c_k(3) = 1$. The slightly less familiar criterion for division into *11* corresponds to $c_k(11) = (-1)^k$.

The universal divisibility rule derived in this work applies to *any* divisor *n* in *any* number system and is of the form

$$C = \sum_{k=0}^{m} u^k w^{m-k} a_k,  \qquad (I.3)$$

where in the *t*-base number system, the parameters *u* and *w* are determined in terms of the divisor *n* or its arbitrary multiple $N = qn$ with integer *q* as follows:

$$N = wt - u. \qquad (I.4)$$

Divisibility parameters *u* and *w* can be found for any *N*, because expression (I.4) is the general form of an arbitrary natural number. They are easily determined for a given *N* (specified by the multiple *q*, positive or negative) but for different *q* one gets different parameters. Thus, different divisibility criteria exist for the same divisor and their choice may be guided by the magnitude of *C* they lead to. The smaller the better. We shall call (I.3) the general divisibility criterion or GDC.

Our method allows one to derive all known divisibility criteria in any base number system.

## II. The general approach; restricted divisibility rules

In this section, we illustrate the general approach by deriving an auxiliary restricted criterion that will be subsequently used in Sect. III to derive the GDC. We shall dispense with superfluous generality and restrict our consideration to the decimal system. Examples will be given in the decimal ($t=10$) and octal ($t=8$) number systems.

Let us represent the test number *A* in a form with separated number of decades *B* and units *b*:

$$A = 10B + b. \qquad (II.1)$$

The restricted divisibility criterion *R* will be constructed as a linear combination of *B* and *b*, viz.

$$R = uB + wb, \qquad (II.2)$$

where *u* and *w* are integer parameters.

Assume that we can choose the divisibility parameters *u* and *w* in such a way that *R* is divisible by *n*. Under this assumption, we shall show that *u* and *w* do not depend on *A* and that the test number *A* itself will be divisible by *n*.

Inasmuch as *R* is divisible by *n*, we can write

$$R = uB + wb = pn,  \qquad (II.3)$$

where *p* is a natural number.

Multiplying (II.1) by *u* and substituting an expression for *uB* from (II.3), we get

$$uA = 10pn - b(10w - u). \qquad (II.4)$$

The first term on the right-hand side of (II.4) is divisible by *n*. Divisibility of the second term can be ascertained by an appropriate choice of parameters *u* and *w*. In other words, the number

$$N = 10w - u \equiv qn \qquad (II.5)$$

must be divisible by *n*. We can take *N* equal to *n* or to an arbitrary multiple of *n*, positive or negative.[2]

Divisibility parameters *u* and *w* can be found for any *N*, because expression $10w - u$ describes an arbitrary natural number (of course, *u* must satisfy the inequality $|u| \leq 9$). For a given divisor, these parameters depend on our choice of *q* and this reflects the possibility of many division criteria of the form (II.2).

Conversely, if *R* given by (II.2) is not divisible by *n* for some choice of *N*, then the test number *A* will not be divisible by *n*.

Let us illustrate the above procedure for *n = 17*. For three choices of *N* (equal to *17, 34,* and *51*) we get, respectively, three pairs of parameters $(w,u)_q$ and three divisibility criteria $R_q$:

$$\begin{aligned}(w,u)_1 &= (2,3); & R_1 &= 3B + 2b \\ (w,u)_2 &= (3,-4); & R_2 &= -4B + 3b \\ (w,u)_3 &= (5,-1); & R_3 &= -B + 5b \end{aligned} \qquad (II.6)$$

At a first glance it may appear that the best choice of parameters *u* and *w* should correspond to *N = n*. However, this is not always true. Even though all three rules (II.6) are correct, the third is the best. The quality of a division criterion is determined by its magnitude, $|R_q|$, the smaller the better. For an exemplary test number A = *5916*, we have $R_1 = 1785$, $R_2 = -2346$, and $R_3 = -561$.

---

[2] Note that multiplication by *u* in Eq. (II.4) would not be meaningful if *u = 0*. However, as seen from (II.5) and (I.4), in this case the number *N* is a multiple of the base *t* of the number system and no divisibility criterion is needed.

Another example corresponds to division by 3. Taking *N = 3, 6,* and *9,* we get

$$(w,u)_1 = (1,7); \quad R_1 = 7B+b$$
$$(w,u)_2 = (1,4); \quad R_2 = 4B+b \quad \quad \quad (\text{II}.7)$$
$$(w,u)_3 = (1,1); \quad R_3 = B+b$$

Again, the third rule (*q* = 3) is best. Generally, the strongest inequality $|R| < A$ obtains for the smallest value of |*u*|, as can be formally shown from Eqs. (II.1) and (II.2). Examples (II.6) and (II.7) illustrate this general rule.

Using the described method, we have created two tables of restricted rules, one for the decimal, the other for the octal system. The tables were generated with minimal effort and are shown here for illustration and reference.

The decimal Table 1 corresponds to the divisor *n* spanned in the interval [3-33] with several numbers added outside this interval. Many of the obtained restricted rules are well-known.

**Table 1. Restricted rules for the decimal number system ( *t = 10* )**

| n | N | u | w | R | n | N | u | w | R | n | N | u | w | R |
|---|---|---|---|---|---|---|---|---|---|---|---|---|---|---|
| 3 | 9 | 1 | 1 | B+b | 16 | -32 | 2 | -3 | 2B-3b | 29 | 29 | 1 | 3 | B+3b |
| 4 | 8 | 2 | 1 | 2B+b | 17 | -51 | 1 | -5 | B-5b | 31 | 31 | -1 | 3 | B-3b |
| 5 | 5 | 5 | 1 | 5B+b | 18 | 18 | 2 | 2 | 2B+2b | 32 | 32 | -2 | 3 | 2B-3b |
| 6 | 12 | -2 | 1 | 2B-b | 19 | 19 | 1 | 2 | B+2b | 33 | -33 | 3 | -3 | 3B-3b |
| 7 | -21 | 1 | -2 | B-2b | 21 | -21 | 1 | -2 | B-2b | 39 | 39 | 1 | 4 | B+4b |
| 8 | 8 | 2 | 1 | 2B+b | 22 | -22 | 2 | -2 | 2B-2b | 49 | 49 | 1 | 5 | B+5b |
| 9 | 9 | 1 | 1 | B+b | 23 | 69 | 1 | 7 | B+7b | 59 | 59 | 1 | 6 | B+6b |
| 11 | 11 | 1 | -1 | B-b | 24 | 48 | 2 | 5 | 2B+5b | 69 | 69 | 1 | 7 | B+7b |
| 12 | -12 | 2 | -1 | 2B-b | 25 | -25 | 5 | -2 | 5B-2b | 79 | 79 | 1 | 8 | B+8b |
| 13 | 39 | 1 | 4 | B+4b | 26 | -52 | 2 | -5 | 2B-5b | 81 | 81 | -1 | 8 | B-8b |
| 14 | 28 | 2 | 3 | 2B+3b | 27 | 27 | 3 | 3 | 3B+3b | 83 | 83 | -3 | 8 | 3B-8b |
| 15 | -15 | 5 | -1 | 5B-b | 28 | 28 | 2 | 3 | 2B+3b | 87 | 87 | 3 | 9 | 3B+9b |

The octal Table 2 comprises divisors *n* spanned in the interval [3-32]$_8$. For several divisors, there are no criteria. This situation occurs for the "round" numbers in the *t = 8* system, such as *$10_8$ = 8, $20_8$ = 16,* and *$30_8$ = 24*. Similar divisors (*10, 20,* and *30*) were omitted in Table 1, see footnote[2].

**Table 2. Restricted rules for the octal number system ( *t = 8* )**

| $n_8$ | $N_8$ | u | w | R | $n_8$ | $N_8$ | u | w | R |
|---|---|---|---|---|---|---|---|---|---|
| 3 | 11 | -1 | 1 | B-b | 17 | 17 | 1 | 2 | B+2b |
| 4 | 10 | 0 | 1 | b | 20 | 20 | 0 | 2 | |
| 5 | 5 | 3 | 1 | 3B+b | 21 | 21 | -1 | 2 | B-2b |
| 6 | 6 | 2 | 1 | 2B+b | 22 | 22 | -2 | 2 | 2B-2b |
| 7 | 7 | 1 | 1 | B+b | 23 | 46 | 2 | 5 | 2B+5b |
| 10 | 10 | 0 | 1 | | 24 | 24 | -4 | 2 | 4B-2b |
| 11 | 11 | -1 | 1 | B-b | 25 | 25 | 3 | 3 | 3B+3b |
| 12 | 12 | -2 | 1 | 2B-b | 26 | 26 | 2 | 3 | 2B+3b |
| 13 | 13 | -3 | 1 | 3B-b | 27 | 27 | 1 | 3 | B+3b |
| 14 | 14 | -4 | 1 | 4B-b | 30 | 30 | 0 | 3 | |
| 15 | 32 | 2 | 3 | 2B+3b | 31 | 31 | -1 | 3 | B-3b |
| 16 | 16 | 2 | 2 | 2B+2b | 32 | 32 | -2 | 3 | 2B-3b |

It is interesting to compare the simplest rules in different base systems. Take the criteria for *n = $3_{10}$, $9_{10}$,* and *$11_{10}$,* cf. Table 1. The first two correspond to division parameters *w = u = 1*, while the third to *w = –u = 1*. The same parameters (cf. Table 2) in the octal system correspond to *$7_8$* and the pair *$3_8$* and *$11_8$* . This means that if in the decimal system the rules for divisors *$3_{10}$* and *$9_{10}$* correspond to divisibility by these numbers of the sum of all digits of the test number, then in the octal system the same rule corresponds to the divisor *$7_8$* . On the other hand, the decimal divisibility criterion for *$11_{10}$* corresponds to divisibility by *$3_8$* and *$11_8$* in the octal system. As if the systems swap their rules!

*Example.* Consider *A = $17223_{10}$*. The number is divisible by *3* because the sum of its digits equals *$15_{10}$*. The same number in the octal system reads *A = $41515_8$*. It is divisible by *3* because the *alternate* sum, *4 – 1 + 5 – 1 + 5 = $14_8$ = $12_{10}$*, of its digits is divisible by *3*.

### III. The general divisibility criterion

Applying the restricted rule described in Sect. II to a test number *A* generates a number *$R_1$* that is divisible by *n* if and only if *A* is divisible by *n*. It is therefore possible to repeat this procedure by applying the restricted rule to the number $R_1$ – thus generating $R_2$ – and so on. This is the usual approach to investigating the divisibility of large numbers with restricted criteria.

We apply this iterative approach to a general (*m+1*)-digit test number *A* (Eq. I.1), taken – without a loss of generality – as a decimal number,

$$A = \sum_{k=0}^{m} 10^k a_k = \sum_{k=1}^{m} 10^k a_k + a_0. \qquad (\text{III.1})$$

First, we re-write Eq. (III.1) in the form (II.1),

$$A = 10 \sum_{k=1}^{m} 10^{k-1} a_k + a_0, \qquad (\text{III.2})$$

and then apply the divisibility rule (II.2) to the number (III.2). We get

$$\begin{aligned} R_1 &= u \sum_{k=1}^{m} 10^{k-1} a_k + w a_0 \\ &= 10 u \sum_{k=2}^{m} 10^{k-2} a_k + (u a_1 + w a_0) \end{aligned} \qquad (\text{III.3})$$

We may call $R_1$ the first-order restricted criterion. Applying the same rule (II.2) to (III.3), we arrive at the second-order restricted criterion,

$$\begin{aligned} R_2 &= u^2 \sum_{k=2}^{m} 10^{k-2} a_k + w(u a_1 + w a_0) \\ &= 10 u^2 \sum_{k=3}^{m} 10^{k-3} a_k + (u^2 a_2 + u w a_1 + w^2 a_0) \end{aligned} \qquad (\text{III.4})$$

This procedure can be continued until all digits of the test number *A* are exhausted. After *m* iterations, we arrive at the following series:

$$C \equiv R_m = \sum_{k=0}^{m} u^{m-k} w^k a_{m-k} = \sum_{k=0}^{m} u^k w^{m-k} a_k, \qquad (\text{III.5})$$

which is a linear form in all digits of the test number *A*. We shall refer to *C* as the general divisibility criterion or GDC.

Note that the form (III.5) of the GDC does not depend on the base *t* of the number system, cf. Eq. (I.3). The division parameters *u* and *w* do, of course, depend on *t* [cf. Eqs. (I.4) and (II.5)].

The procedure we just used to derive Eq. (III.5) requires a clarification. The validity of applying rule (II.2) to the algebraic expression of restricted rules $R_k$ is in need of proof. Expression (II.2) uses a linear combination of the decades and the units of the test number. There is no objection to applying (II.2) to *A*.

However, we cannot *a priori* separate the number $R_1$ into decades and units according to the decimal rules. We have represented $R_1$ in such a form only formally – without checking that the term $(ua_1 + wa_0)$ is smaller than *9*. In general, it is not. The same problem arises in higher orders of $R_k$. No such problem arises in numerical calculations, where ordering of the decimal numbers is automatically regulated by the rules of arithmetic calculations.

To show that the invalid form of $R_k$ does not destroy the divisibility of these numbers, we consider the test number *A* in the form (II.1) and suppose that the number of units *b* contains *K* extra decades. To "correct" this invalid decimal representation, we subtract the extra decades from *b* and add them to *B*, viz.

$$A = 10(B+K) + (b-10K).  \qquad (III.6)$$

Applying the divisibility criterion (II.2) to Eq. (III.6), we obtain

$$\begin{aligned} R &= u(B+K) + w(b-10K) \\ &= uB + wb - K(10w - u) \end{aligned}. \qquad (III.7)$$

The extra term in (III.7) equals *KN*, where *N* is defined by (II.5) and is a multiple of the divisor *n*. This proves that the invalid decimal representation of *R* does not destroy its divisibility.

**IV. Discussion and conclusion**

The GDC acquires its simplest form, when the parameters *u* and *w* have equal absolute values. This is true for any number system. Equations (I.3) and (I.4) yield

$$\begin{aligned} C &= |w|^m \sum_{k=0}^{m} \left( \frac{uw}{|uw|} \right)^k a_k \\ N &= (t \pm 1)w \end{aligned}. \qquad (IV.1)$$

Thus, in the decimal system for multiples of 3 and 9, both divisibility parameters can be taken equal,[3] *u = w*, and this yields

---

[3] We remark that the divisibility parameters (*w,u*) = (*8,--1*) chosen in Table 1 for *n = 81*, were optimized (taking the set with smallest |*u*|) for the *restricted* criterion. For the GDC one could advantageously choose the two-parameter set *u = w = 9* and then use (IV.2). It is worthwhile to point out that both parameter sets (*w,u*) = (*8,--1*) and (*w,u*) = (9, 9) correspond to the same divisor multiple *N = 81*, i.e. same *q = 1*, cf. Eq. (II.5).

Note that two different sets of divisibility parameters can be generated using the same *N*. This is due to the fact that expression (I.4) allows both positive and negative values of *u* and hence two representations of the number *N*, e.g., 81 = 8·10 + 1 = 9·10 – 9.

$$C = u^m \sum_{k=0}^{m} a_k . \tag{IV.2}$$

On the other hand, for multiples of *11*, one can take *u = - w* to obtain

$$C = u^m \sum_{k=0}^{m} (-1)^k a_k . \tag{IV.3}$$

The same expressions (IV.2) and (IV.3) are valid in the octal system (*t = 8*) but, respectively, for multiples of $7_8$ and $11_8$. This "reciprocity" of rules was already noted in Sect. II underneath Table 2.

Another simplification of the GDC arises when the test number *A* has identical digits, $a_k = a$. In this case, the rule (III.5) reduces to

$$C = a \frac{u^{m+1} - w^{m+1}}{u - w} = a \sum_{k=0}^{m} u^k w^{m-k} . \tag{IV.4}$$

A curious reader may be tempted to derive from the GDC multiple special rules for selected divisors, test numbers and base systems.

The GDC power can be demonstrated by considering an example of divisor *7* in the decimal system. Parameters *u* and *w* in this case are (see Table 1): *u = 1* and *w = –2*. For an exemplary test number *A = 1860523*, Eq. (III.5) yields *C = 217*. This number is about four orders of magnitude smaller than *A*.

In conclusion, we have demonstrated a general divisibility criterion (GDC) of integers. It is universal in the sense that it can be applied to any divisor in any base number system. The criterion is written down as a linear form in all digits of the test number and is determined by two parameters characterizing the divisor in a given number system.

The direct usefulness of the GDC for testing the divisibility of numbers may not be that great in the age of computers. It would perhaps be more appreciated for this purpose by ancient Greeks... However, due to its explicit nature, the GDC may become useful as an analytical tool in a variety of mathematical problems.

## A. Appendix: Simplified derivation of the GDC

This note is added in March 2016. Using the notation of congruence arithmetic, the general divisibility criterion derived in this paper can be written in the form

$$A = \sum_{k=0}^{m} t^k a_k \equiv \sum_{k=0}^{m} u^k w^{m-k} a_k \; [\mathrm{mod}\, N] = \mathrm{GDC}, \tag{A.1}$$

where, in the *t*-base number system, the divisor *N* is given by

$$N = wt - u. \tag{A.2}$$

To prove (A.1), we note the congruence relation,

$$(wt)^k = (wt - u + u)^k \equiv (u)^k \; [\mathrm{mod}\,(wt - u)]. \tag{A.3}$$

Next, we divide and multiply the left-hand side of (A.1) by $w^m$ and then use the relation (A.3). This gives

$$\begin{aligned}
A &= \frac{1}{w^m} \sum_{k=0}^{m} (wt)^k w^{m-k} a_k \\
&\equiv \frac{1}{w^m} \sum_{k=0}^{m} u^k w^{m-k} a_k \quad [\mathrm{mod}\, N] \\
&\equiv \sum_{k=0}^{m} u^k w^{m-k} a_k \quad [\mathrm{mod}\, N] \\
&= \mathrm{GDC}
\end{aligned} \tag{A.4}$$

This means that *A* is divisible into *N* if and only if the GDC is divisible into *N*.